\documentclass{article}

\usepackage[alphabetic]{amsrefs}
\usepackage{amsmath,amssymb,enumerate}
\usepackage[all,cmtip]{xy}

\def\today{\number\day .\number\month .\number\year}
\usepackage[applemac]{inputenc}

\def\1{{\bf 1}}
\def\A{{\mathbb A}}
\def\Ad{\operatorname{Ad}}
\def\al{\alpha}
\def\bs{\backslash}
\def\C{{\mathbb C}}
\def\CA{{\cal A}}
\def\CC{{\cal C}}
\def\CD{{\cal D}}
\def\CF{{\cal F}}
\def\CG{{\cal G}}
\def\CJ{{\cal J}}
\def\CO{{\cal O}}
\def\CP{{\cal P}}

\def\cusp{{\rm cusp}}
\def\comm{\operatorname{comm}}
\def\df{\ \stackrel{\mbox{\rm\tiny def}}{=}\ }
\def \ds{\displaystyle}
\def\e{\emph}
\def\eps{\varepsilon}
\def\Ext{\operatorname{Ext}}
\def\fin{{\rm fin}}
\def\g{{\mathfrak g}}
\def\G{{\mathbb G}}
\def\ga{\gamma}
\def\Ga{\Gamma}
\def\GL{\operatorname{GL}}
\def\Gr{{\rm Gr}}
\def\H{\operatorname{H}}
\def\HH{{\mathbb H}}
\def\Harm{\operatorname{Harm}}
\def\Hom{{\rm Hom}}
\def\Im{\operatorname{Im}}
\def\k{{\mathfrak k}}
\def\La{\Lambda}
\def\la{\lambda}
\def\mathqed{\vspace{-30pt}\\ \qed}
\def\mg{{\rm mg}}
\def\Mod{\operatorname{Mod}}
\def\N{{\mathbb N}}

\def\p{{\mathfrak p}}

\def\pa{{\rm par}}
\def\Pet{{\rm Pet}}
\def\PGL{\operatorname{PGL}}
\def\ph{\varphi}
\def\prf{{\bf Proof: }}
\def\PSL{\operatorname{PSL}}
\def\Q{{\mathbb Q}}
\def\qed{\ifmmode\eqno \square 
		\else\noproof\vskip 12pt plus 3pt minus 9pt \fi}
\def\noproof{{\unskip\nobreak\hfill\penalty50\hskip2em\hbox{}%
     \nobreak\hfill $\square$\parfillskip=0pt%
     \finalhyphendemerits=0\par}}
\def\R{{\mathbb R}}

\def\res{\operatorname{res}}
\def\SL{{\rm SL}}
\def\sm{\smallsetminus}
\def\triv{{\rm triv}}
\def\ul{\underline}
\def\umg{{\rm umg}}
\def\vol{\operatorname{vol}}
\def\what{\widehat}
\def\Z{{\mathbb Z}}
\def\z{{\mathfrak z}}
\def\({\left(}
\def\){\right)}

\newcommand{\tto}[1]{\stackrel{#1}{\longrightarrow}}

\newcommand{\norm}[1]{|\hspace{-1pt}| #1|\hspace{-1pt}|}
\renewcommand{\sp}[1]{\left\langle #1\right\rangle}
\newcommand{\ol}[1]{\overline{#1}}
\newcommand{\stack}[2]{\genfrac{}{}{0pt}{1}{#1}{#2}}

\newcommand{\smat}[4]{\(\begin{smallmatrix}#1 & #2 \\ #3 & #4\end{smallmatrix}\)}

\newtheorem{theorem}{Theorem}[subsection]
\newtheorem{conjecture}[theorem]{Conjecture}
\newtheorem{lemma}[theorem]{Lemma}

\newtheorem{proposition}[theorem]{Proposition}

\begin{document}

\pagestyle{myheadings} \markright{INVARIANTS, COHOMOLOGY, AND AUTOMORPHIC FORMS...}

\title{Invariants, cohomology, and  automorphic forms of higher order}
\author{Anton Deitmar\\ \ \\ 
Selecta Math. New Ser. 18, No 4, 855-883 (2012)}
\date{}
\maketitle

{\bf Abstract:}
A general structure theorem on higher order invariants is proven.
For an arithmetic group, the structure of the corresponding Hecke module is determined.
It is shown that the module does not contain any irreducible submodule.
This explains the fact that L-functions of higher order forms have no Euler-product.
Higher order cohomology is introduced, classical results of Borel are generalized and a higher order version of Borel's conjecture is stated.




\section*{Introduction}

Higher order modular forms show up in various contexts.
For instance, the limit crossing probabilities in percolation theory, considered as functions of the aspect ratio, turn out to be higher order forms \cite{KZ}.
As an approach to the abc-conjecture, Dorian Goldfeld introduced Eisenstein-series twisted by modular symbols \cite{CDO,GG,Gold}, which are higher order forms.
Finally, spaces of higher order forms are natural receptacles of converse theorems \cite{F,IM}.

L-functions of second order forms have been studied in \cite{DKMO}, Poincar\'e series attached to higher order forms have been investigated in \cite{IO}, dimensions of spaces of second order forms have been determined in \cite{DO,DiamSim}.
Higher order cohomology has been introduced and an Eichler-Shimura type theorem has been proven in \cite{ES}.
In \cite{HOAut} a program has been started which aims at an understanding of the theory of higher order forms from an representation-theoretical point of view.
The present paper contains a major step in that direction.
It is shown that the order-lowering-homomorphism introduced in \cite{HOAut} injects the graded module of higher order forms into a canonical tensor product.
In the case of cusp forms of weight two and order one it is shown that the order-lowering-homomorphism is indeed an isomorphism.
So in that case the representation theoretic nature of higher order forms is completely understood.
In the general situation it is not clear under which circumstances the order-lowering-map is an isomorphism.

In this paper we also resolve a standing issue of higher order forms: we give a natural action of the Hecke algebra.
Using the previous results, we find that this action has no eigenvectors, which explains the absence of Euler-products for $L$-functions of higher order forms.
The Hecke action is derived from an action of the group $\CG=G(\Q)$ of rational points of the ambient reductive linear group $G$.
This action does not preserve a given arithmetic group, so we resort to taking limits over all congruence subgroups which naturally leads to adelic groups and their representations.
We obtain a map that injects the space of higher order forms into the tensor product of automorphic representations.
More precisely, if $S_{k,q}$ is the space of cusp forms of weight $k$ and order $q$, we consider the graded piece 
$$
\Gr S_{k,q}=S_{k,q}/S_{k,q-1}
$$
and obtain an map
$$
\psi : \Gr S_{k,q}\hookrightarrow (S_1\oplus \ol{S_2})^{\otimes (q-1)}\otimes S_k.
$$
We show that there are natural pre-Hilbert space structures on these spaces making $\psi$ an isometry.
In the case $k=2=q$ we show in Theorem \ref{thm1.4.2} that, after completion, $\psi$ is an isomorphism, which is striking, as before taking the limit over all congruence subgroups,  $\psi$ is never surjective.
It is an open question whether $\psi$ is an isomorphism for other values of $k$ and $q$.
We next show that a tensor product of two automorphic representations has no irreducible subrepresentation and conclude that $S_{k,q}$ has no irreducible subrepresentation.
This explains the empirical fact that L-functions of higher order forms do not admit Euler products.

In the second part of the paper we relate higher order forms to higher order cohomology of arithmetic groups.
It turns out that on the level of Lie algebra cohomology it is possible to separate the steps of taking higher order invariants from the step of taking cohomology.
WE give a full spectral decomposition of the space of higher order forms in the cocompact case and we formulate the higher order version of the Borel conjecture which asks whether higher order cohomology of arithmetic groups can be computed from the complex of higher order automorphic forms.

\section{Higher order invariants}
\subsection{Generalities}
Let $R$ be a commutative ring with unit,
let $\Ga$ be a group and let $P\subset \Ga$ be a conjugation-invariant subset.
Let $I$ denote the augmentation ideal in the group algebra $A=R[\Ga]$.
For an $A$-module $V$ and $q\in\N$ we define the $R$-module
$$
\H_{q}^0(\Ga,P,V)
$$
of invariants of order $q$ to be the set of all $v\in V$ with $I^{q}v=0$ and $pv=v$ for every $p\in P$.
Note that for $q=1$ one gets the usual invariants
$$
\H_{1}^0(\Ga,P,V)= \H^0(\Ga,V)= V^\Ga.
$$
Another way to describe the space $\H_{q}^0(\Ga,P,V)$ is as follows.
Let $\sp P\subset \Ga$ be the subgroup generated by $P$ and let $I_P$ be the augmentation ideal of the group $\sp P$.
As $P$ is conjugation-invariant, the group $\sp P$ is normal in $\Ga$ and so $AI_P$ is a two-sided ideal in $A$.
Let
$$
J_{q}\df I^{q}+AI_P.
$$
Then $J_{q}$ is a two-sided ideal of $A$ and $\H_{q}^0(\Ga,P,V)$ is the space of all $v\in V$ with $J_{q}v=0$.
There is a natural identification
$$
\H_{q}^0(\Ga,P,V)\ \cong\ \Hom_A(A/J_{q},V).
$$
If $P$ is the empty set or $\{ 1\}$, we write $\H_q^0(\Ga,V)$ for $\H_{q}^0(\Ga,P,V)$.

The sets $\H_{q}^0(\Ga,P,V)\subset 
\H_{q+1}^0(\Ga,P,V)$ form a filtration on $V$ which is not necessarily exhaustive.
Let
$$
\ol\H_{q}(\Ga,P,V)\df \H_{q}^0(\Ga,P,V)/\H_{q-1}^0(\Ga,P,V)
$$
be the $q$-th graded piece, where we allow $q=1,2,\dots$ by formally setting $\H_{0}^0(\Ga,P,V)=0$.

For any group $G$, let $\Hom_P(\Ga,G)=\Hom(\Ga/\sp P,G)$.
We will in particular use this notation for $G$ being the additive group of a vector space $W$.

If $R$ is a field and if $\Ga$ is finitely generated or $\dim W<\infty$, then
$$
\Hom_P(\Ga,W)\cong \Hom_P(\Ga,\C)\otimes W,
$$
where the tensor product is over $R$.
We will make use of this fact in the sequel.
We introduce the \e{order-lowering-homomorphism}
$$
\La: \ol\H_{q}(\Ga,P,V)\to\Hom_P(\Ga,\ol\H_{q-1}(\Ga,P,V))
$$
given by
$$
\La(v)(\ga)=(\ga-1)v.
$$
To see that $\La(v)$ is indeed a group homomorphism, note that in the group algebra $\C[\Ga]$ one has $(\tau\ga-1)\equiv (\tau-1)+(\ga-1)\mod I^2$.

\begin{lemma}
Let $V$ be an $R[\Ga]$-module which is torsion-free as $\Z$-module.
Let $\Sigma\subset\Ga$ be a subgroup of finite index.
Then the natural restriction map
$$
\res: \ol\H_{q}(\Ga,P,V)\ \to\ \ol \H_{q}(\Sigma,P\cap\Sigma,V)
$$
is injective.
\end{lemma}

\prf
For convenience, we will also write $\H_{q}(\Sigma,P,V)$ instead of $\H_{q}(\Sigma,P\cap\Sigma,V)$.
We prove the lemma by induction on $q$.
For $q=1$ the restriction map is the inclusion $V^\Ga\hookrightarrow V^\Sigma$ and thus injective.
For $q\ge 2$ we have the commutative diagram
$$
\xymatrix{
\ol\H_{q}(\Ga,P,V)\ar[r]^{\res} \ar@{^{(}->}[d]_\La
& \ol \H_{q}(\Sigma,P,V)\ar@{^{(}->}[d]_\La
\\
\Hom_P(\Ga,\ol \H_{q-1}(\Ga,P,V))\ar[r]^\res
& \Hom_P(\Sigma,\bar H_{q-1}(\Sigma,P,V)).
}
$$
To see that the top row is injective, we have to show that the bottom row is injective.
The bottom row is the composition of two maps
\begin{eqnarray*}
\Hom_P(\Ga,\ol \H_{q-1}(\Ga,P,V))&\to&
\Hom_P(\Sigma,\ol \H_{q-1}(\Ga,P,V))\\
&\to&
\Hom_P(\Sigma,\ol \H_{q-1}(\Sigma,P,V)),
\end{eqnarray*}

the first of which is injective as $\Sigma$ is of finite index and $V$ is torsion-free, and the second is injective by induction hypothesis.
\qed

\subsection{Hecke pairs and smooth modules}
A \e{Hecke pair} is a pair $(\CG,\Ga)$ of a group $\CG$ and a subgroup $\Ga$ such that for every $g\in \CG$ the set $\Ga g\Ga/\Ga$ is finite.
We also say that $\Ga$ is a \e{Hecke subgroup} of $\CG$.

Two subgroups $\Ga,\La$ of a group $H$ are called \e{commensurable}, written $\Ga\sim\La$,  if the intersection $\Ga\cap\La$ has finite index in both.
Commensurability is an equivalence relation which is preserved by automorphisms of $H$.

The \e{commensurator} of a group $\Ga\subset H$ is
$$
\comm(\Ga)\df \{ h\in H:\Ga \text{ and } h\Ga h^{-1}\text{ are commensurable}\}.
$$
The commensurator $G=\comm(\Ga)$ is a subgroup of $H$.
It is the largest subgroup such that $(G,\Ga)$ is a Hecke pair.
More precisely,
$$
\comm(\Ga)=\{ h\in H: |\Ga h\Ga/\Ga|,|\Ga\bs\Ga h\Ga|<\infty\}.
$$

Let $\CG$ be a group.
By a $\CG$-module we shall henceforth mean an $R[\CG]$-module.
If $\CG$ is a totally disconnected topological group, an element $v$ of a $\CG$-module $V$ is called \e{smooth} if it is stabilized by some open subgroup of the topological group $\CG$.
The set $V^\infty$ of all smooth elements is a submodule and the module $V$ is called \e{smooth} if $V=V^\infty$.

Drop the condition that $\CG$ be a topological group and let $(\CG,\Ga)$ be a Hecke-pair.
A \e{congruence subgroup} of $\Ga$ is any subgroup which contains a group of the form
$$
\Ga\cap g_1\Ga g_1^{-1}\cap\dots\cap g_n\Ga g_n^{-1}
$$
for some $g_1,\dots, g_n\in \CG$.
As $(\CG,\Ga)$ is a Hecke pair, every congruence subgroup has finite index in $\Ga$.
Note that this definition of a congruence subgroup coincides with the one given in \cite{HOAut}.
For every congruence subgroup $\Sigma$ equip the set $\CG/\Sigma$ with the discrete topology and consider the topological space
$$
\bar\CG\df \lim_{\stack\leftarrow\Sigma}\CG/\Sigma,
$$
where the limit is taken over all congruence subgroups $\Sigma$.

\begin{lemma}
\begin{enumerate}[\rm (a)]
\item The intersection of all congruence subgroups $N=\bigcap_\Sigma\Sigma$ is a normal subgroup of $\CG$.
\item The natural map $p:\CG\to\bar\CG$ factors through the injection $\CG/N\hookrightarrow \bar \CG$ and has dense image.
\item The group multiplication on $\CG/N$ extends by continuity to $\bar \CG$ and makes $\bar\CG$ to a totally disconnected locally compact group.
\end{enumerate}
\end{lemma}

We call $\bar\CG$ the \e{congruence completion} of $\CG$.
Although the notation doesn't reflect this, the completion $\bar\CG$ depends on the choice of the Hecke subgroup $\Ga$.
A Hecke subgroup $\Ga$ is called \e{effective}, if the normal subgroup $N$ above is trivial.

\prf
(a) Let $n\in N$ and let $g\in \CG$.
For a given congruence subgroup $\Sigma$ we have that $n\in\Sigma\cap g^{-1}\Sigma g$, and so $gng^{-1}\in\Sigma$.
As $\Sigma$ varies, we find $gng^{-1}\in N$.

(b) Let $g,g'\in\CG$ with $p(g)=p(g')$.
This means that $g\Sigma=g'\Sigma$ for every congruence subgroup and so $gN=g'N$. 
For given $(g_\Sigma)_\Sigma\in\bar\CG$ the sets $U_\Sigma=\{ h\in\bar\CG: h_\Sigma\Sigma=g_\Sigma\Sigma\}$ form a neighborhood base.
Clearly the element $g_\Sigma\in\CG$ is mapped into $U_\Sigma$, so the image of $p$ is dense.

(c) Let $\bar g=(g_\Sigma)_\Sigma\in\bar\CG$.
Then the net $(p(g_\Sigma))_\Sigma$ converges to $\bar g$.
For $\bar h=(h_\Sigma)\in\bar\CG$ it is easy to see that the net $(p(g_\Sigma h_\Sigma))_\Sigma$ converges in $\bar\CG$.
We call the limit $\bar g\bar h$.
This multiplication has the desired properties.
\qed

The initial topology defined by $p$ on $\CG$ makes $\CG$ a topological group with the congruence groups forming a unit neighborhood base.
Clearly every smooth $\bar\CG$-module is a smooth $\CG$-module by restriction.
But also the converse is true: Every smooth $\CG$-module extends uniquely to a smooth $\bar\CG$-module and these two operations of restriction and extension are inverse to each other.

Let $\CP\subset\CG$ be a conjugation-invariant set.
For any subgroup $\Sigma$ of $\Ga$ we write $\H_{q}^0(\Sigma,\CP,V)$ for $\H_{q}(\Sigma,\CP\cap\Sigma,V)$.
Let
$$
L_{q}(\CP,V)\df \lim_{\stack\to\Sigma}\ol \H_{q}(\Sigma,\CP,V),
$$
where the limit is taken over all congruence subgroups of $\Ga$.
Note that $L_{1}(\CP,V)=V^\infty$.
For $g\in\CG$, the map induced by $g$:
$$
\ol \H_{q}(\Sigma,\CP,V)\to \ol \H_{q}^0(g\Sigma g^{-1},\CP,V)\tto\res\ol \H_{q}^0(\Sigma\cap g\Sigma g^{-1},\CP,V)
$$
defines an action of $\CG$ on $L_{q}(\CP,V)$, which makes the latter a smooth module.

Assume from now on that $\Ga$ is finitely generated and $R$ is a field.
Then every finite-index subgroup $\Sigma\subset\Ga$ is finitely generated as well.

Consider the order-lowering map
$$
\ol \H_{q}(\Sigma,\CP,V)\hookrightarrow\Hom_\CP(\Sigma,\C)\otimes\ol \H_{q-1}(\Sigma,\CP,V),
$$
where $\Hom_\CP(\Sigma,\C)$ is the set of all homomorphisms $\Sigma\to\C$ that vanish on $\Sigma\cap \CP$.
Iteration gives
$$
\ol \H_{q}(\Sigma,\CP,V)\hookrightarrow\Hom_\CP(\Sigma,\C)^{\otimes (q-1)}\otimes V^\Sigma.
$$
Taking limits we get an injection
$$
L_{q}(\CP,V)\hookrightarrow \(\lim_{\stack\to\Sigma}\Hom_\CP(\Sigma,\C)\)^{\otimes (q-1)}\otimes V^\infty.
$$
Let $\what\Hom_{\CP}$ be the space $\displaystyle\lim_{\stack\to\Sigma}\Hom_\CP(\Sigma,\C)$.
For $g\in\CG$ one gets a map
$$
\Hom_\CP(\Sigma,\C)\to\Hom_\CP(g\Sigma g^{-1},\C)\tto\res\Hom_\CP(\Sigma\cap g\Sigma g^{-1},\C),
$$
which makes $\what\Hom_{\CP}$ a smooth module.
We have shown:

\begin{proposition}\label{4.2}
If $\Ga$ is finitely generated and $R$ is a field,
then there is a natural injection of smooth modules 
$$
L_{q}(\CP,V)\ \hookrightarrow\ 
\what\Hom_{\CP}^{\otimes (q-1)}\otimes V^\infty.
$$ 
\end{proposition}

\subsection{Arithmetic groups}
The ring $R$ continues to be the field $\C$.
Let $\G$ be a linear algebraic group over $\Q$ which is simple and simply connected and such that $G=\G(\R)$ has no compact component.
By strong approximation, the group $\CG=\G(\Q)$ is dense in $\G(\A_\fin)$.
Let $K_\fin$ be a given compact open subgroup of $\G(\A_\fin)$ and let $\Ga=\CG\cap K_\fin$.
Then $(\CG,\Ga)$ is a Hecke-pair and $\Ga$ can be chosen to be effective.
In this case the congruence completion is
$$
\bar\CG\ \cong\ \G(\A_\fin).
$$
We want to analyze the smooth module 
$$
\what\Hom_\Ga=\lim_{\stack\to\Sigma}\Hom(\Sigma,\C).
$$ 
We first consider the case when the $\Q$-rank of $\G$ is zero.
Then $\Ga$ is cocompact in $G$.
The $\G(\A)$-representation on $L^2(\G(\Q)\bs\G(\A))$ decomposes into a direct sum
$$
\bigoplus_{\pi\in\widehat{\G(\A)}}N(\pi)\pi,
$$
where $\widehat{\G(\A)}$ is the unitary dual of $\G(\A)$ and $N(\pi)\in\N_0$ is the multiplicity of $\pi$.
Note that every $\pi\in\widehat{\G(\A)}$ is a tensor product
$$
\pi=\(\bigotimes_{p}\pi_p\)\ \otimes\ \pi_\infty,
$$
where the product runs over all primes $p$ and $\pi_p\in\widehat{\G(\Q_p)}$.
We also denote the representation $\bigotimes_{p}\pi_p$ of the group $\G(\A_\fin)$ by $\pi_\fin$.

Let $\g$ be the complexified Lie algebra of $G$ and $\k$ the complexified Lie algebra of $K$.
Then $K$ acts on $\g/\k$ by the adjoint representation.

\begin{theorem}\label{4.3}
If the $\Q$-rank of $\G$ is zero, then $\what\Hom_\Ga$ is the space of smooth vectors $\(\bar H_\Ga\)^\infty$ of a unitary representation $\bar H_\Ga$, which
decomposes into a direct Hilbert sum of irreducible representations with finite multiplicities,
$$
\bar H_\Ga\ \cong\ \bigoplus_{\pi_\fin\in\widehat{\G(\A_\fin)}}m(\pi_\fin)\ \pi_\fin.
$$
The multiplicity $m(\pi_\fin)$ equals
$$
m(\pi_\fin)= 
\sum_{\stack{\pi_\infty\in\widehat G}{\pi_\infty(C)=0}}d(\pi_\infty) N(\pi_\fin\otimes\pi_\infty),
$$
where $C$ is the Casimir operator of $G$, and $d(\pi_\infty)=\dim\Hom_K(\g/\k,\pi_\infty)$.
The  sum is finite as there are only finitely many $\pi_\infty$ with $\pi_\infty(C)=0$ and $d(\pi_\infty)\ne 0$.
\end{theorem}

\prf
There exists a congruence subgroup $\Sigma\subset\Ga$ which is torsion-free.
For every such $\Sigma$ we have
$$
\Hom(\Sigma,\C)\cong \H^1(\Sigma,\C)\cong \H^1(\Sigma\bs X,\C)\cong \Harm^1(\Sigma\bs X),
$$
where $X= G/K$ is the symmetric space attached to $G$, and $\Harm^1(\Sigma\bs X)$ is the space of harmonic $1$-forms.
The latter space identifies naturally with
$$
W=\{ f\in (C^\infty(\Sigma\bs G)\otimes\g/\k)^K: Cf=0\}.
$$
Now 
$$
C^\infty(\Sigma\bs G)\subset L^2(\Sigma\bs G)=L^2(\G(\Q)\bs\G(\A))^{K_\Sigma}=\bigoplus_{\pi\in\widehat{\G(\A)}}N(\pi)\ \pi^{K_\Sigma}.
$$ 
So that
$$
W= \bigoplus_{\stack{\pi\in\widehat{\G(\A)}}{\pi_\infty(C)=0}}N(\pi)\ (\pi_\infty\otimes\g/\k)^K\otimes\pi_\fin^{K_\Sigma}
$$
The claim follows.
\qed

We do not expect this result to hold for $\Q$-rank$(G)\ge 1$, but in that case we expect a similar result for $\bar H_{\Ga,\CP}$, where $\CP\subset \G(\Q)$ is the union of all $P^1\cap\G(\Q)$, where $P=MAN$ runs over all parabolic subgroups of $G$ which are defined over $\Q$ and $P^1=MN$.
In the case of the group $\G=\PGL_2$ this fact is a consequence of the Eichler-Shimura Theorem.

\subsection{Cusp forms}
We now write $\G$ for the linear algebraic group $\PGL_2$.
Let $\HH$ be the upper half plane in $\C$.
The group $G=\G(\R)$ acts on $\C\sm\R=\HH\cup\ol\HH$.

For $\al=\smat abcd\in\GL_2(\R)$, $k\in 2\Z$, and $f:\C\sm\R\to\C$ we write
$$
f|_k\al(z)=|\det(\al)|^{k/2}(cz+d)^{-k}f\(\frac{az+b}{cz+d}\).
$$
This defines a right action of the group $G$ on the space $\CO(\C\sm\R)$ of holomorphic functions on $\C\sm\R$.
The group $G^+\subset G$ of all matrices of positive determinant is the connected component of $G$.
It is isomorphic to $\PSL_2(\R)$.
The action of $G^+$ leaves stable the upper half plane $\H$, and so the $|_k$-action leaves stable the space $\CO(\H)$.

The space $\CO(\ol\HH)$ can be identified with the space $\ol{\CO(\HH)}$ by the map $f(z)\mapsto f(\ol z)$.
So the group $G$ acts on $\CO(\HH)\oplus\ol{\CO(\HH)}$.
This action can be described by
$$
f|_k\al(z)=
\begin{cases}|\det(\al)|^{k/2}(cz+d)^{-k}f\(\frac{az+b}{cz+d}\)
& \text{if }\det(\al)>0,\\
|\det(\al)|^{k/2}(c\ol z+d)^{-k}f\(\frac{a\ol z+b}{c\ol z+d}\)
& \text{if }\det(\al)<0.\end{cases}
$$

We consider two submodules.
For every cusp $c\in\R\cup\{\infty\}$ of $\Ga=\PSL_2(\Z)$ we fix $\eta_c\in G$ with $\eta_c(\infty)=c$ and $\eta_c^{-1}\Ga_c\eta_c=\pm\smat {1}\Z\ { 1}$.
We define $\CO_{\Ga,k}^M(\H)$ to be the space of all $f\in\CO(\H)$ such that for every cusp $c$ the function $f|_k{\eta_c}(z)$ is, in the domain $\{\Im(z)>1\}$, bounded by a constant times $\Im(z)^A$ for some constant $A>0$.

Further, we define $\CO_{\Ga,k}^S(\HH)$ to be the set of all $f\in\CO(\HH)$ such that for every cusp $c$ the function $f|_{\eta_c}(z)$ is, in the domain $\{\Im(z)>1\}$, bounded by a constant times $e^{-A\Im(z)}$ for some constant $A>0$.
We then have
$$
\CO_{\Ga,k}^S(\HH)\ \subset\ \CO_{\Ga,k}^M(\HH)\ \subset\ \CO(\HH).
$$

We let $\CG=\G(\Q)^+=\G(\Q)\cap G^+$.
For a set $S$ of places of $\Q$, let $\A_S$ be the ring of $S$-adeles, i.e., the ring of all $x\in\prod_{p\in S}\Q_p$ such that $x_p\in\Z_p$ for almost all $p$.

For a ring $R$ we write $\PGL_2(R)$ for the group $\GL_2(R)/R^\times$.
This notion is slightly misleading as the group of $R$-points of the group scheme $\PGL_2$ does not for all rings $R$ coincide with our notion, but one can show that it does so for fields, for $R=\Z_p$ as well as $R=\A_S$ for any set of places $S$.

\begin{lemma}
The group $\G(\Q)^+$ is dense in $\G(\A_\fin)$.
\end{lemma}

\prf
By strong approximation for $\SL_2$ together with  the same assertion for $\GL_1$, it follows that the group $\GL_2(\A)$ is dense in $\GL_2(\A_\fin)$.
As $\GL_2(\Q)^+$ has index two in $\GL_2(\Q)$, its closure in $\GL_2(\A_\fin)$ has at most index 2.
In order to show that the index is one, it suffices to show that the closure contains $T=\smat {-1}\ \ 1$.
For given distinct primes $p_1,\dots,p_k$ and natural numbers $n_1,\dots n_k$ one considers the number $m=p_1^{n_1}\cdots p_k^{n_k}-1$.
Then $m\in\N$ and $m\equiv -1\mod p_j^{n_j}$ for every $j=1,\dots,k$.
For varying $p_j$ and $n_j$ the matrices $\smat m\ \ 1$ will approximate $T$ in $\G(\A_\fin)$.
This shows the claim.
\qed

We define $\CP\subset \CG$ to be the subset of all parabolic elements.
In this case we also use the index ``par'' instead of ``$\CP$''.
Since $\G(\Q)^+\G(\R)$ is dense in $\G(\A)$ it is still true that
$
\bar\CG\ \cong\ \PGL_2(\A_\fin).
$

For a congruence subgroup $\Sigma\subset\Ga$ we define 
$$
M_{k,q}(\Sigma)\df \H_q^0(\Sigma,\CO_{\Ga,k}^M(\HH)).
$$
This is the space of higher order modular forms.
The space of higher order cusp forms is defined as
$$
S_{k,q}(\Sigma)\df H_{q}^0(\Sigma,\Sigma_\pa,\CO_{\Ga,k}^S(\HH)).
$$
The filtration $S_{k,q}(\Sigma)\subset S_{k,q+1}(\Sigma)$ defines the graded pieces
$$
{\rm Gr}S_{k,q}(\Sigma)= S_{k,q}(\Sigma)/S_{k,q-1}(\Sigma).
$$
We also denote 
$$
\Gr S_{k,q}=\lim_{\stack\to\Sigma} \Gr S_{k,q}(\Sigma).
$$

On the space $S_2(\Sigma)=S_{2,1}(\Sigma)$ we introduce the \e{Petersson inner product},
$$
\sp{f,g}_\Pet=\frac 1{[\Ga:\Sigma]}\int_{\Sigma\bs\HH} f(z)\ol{g(z)}\,dx\,dy.
$$
It defines an inner product on the space $S_2=\lim_\Sigma S_2(\Sigma)$.
In the same way, one gets an inner product on the complex conjugate $\ol{S_2(\Sigma)}$ and hence on $\(S_2(\Sigma)\oplus\ol{S_2(\Sigma)}\)\otimes S_2(\Sigma)$.

The Eichler-Shimura isomorphism gives 
$$
\Hom_\pa(\Sigma,\C)\cong S_2(\Sigma)\oplus\ol{S_2(\Sigma)},
$$
 and so $\Gr S_{2,2}(\Sigma)$ injects into $\(S_2(\Sigma)\oplus\ol{S_2(\Sigma)}\)\otimes S_2(\Sigma)$.
Taking limits we get a natural map from $\Gr S_{2,2}$ to $(S_2\oplus \ol{S_2})\otimes S_2$.
This defines an inner product on $\Gr S_{2,2}$.
We denote by $\widehat{\Gr S}_{2,2}$ the Hilbert space completion of $\Gr S_{2,2}$.

The following theorem gives the structure of $\widehat{\Gr S}_{2,2}$ as $\bar\CG$-module and thus as a Hecke-module (compare \cite{HOAut}).
It therefore answers the principal question raised in that paper.

\begin{theorem}\label{thm1.4.2}
The natural map
$$
\psi:\widehat{\Gr S}_{2,2}\ \to\ \( \widehat S_2\oplus \widehat{\ol{S_2}}\)\otimes \widehat S_2
$$
is a unitary $\G(\A_\fin)$-isomorphism.
\end{theorem}

\prf
For a congruence subgroup $\Sigma$ of $\Ga=\PSL_2(\Z)$ let $K_\Sigma$ be the closure of $\Sigma$ in $\G(\A_\fin)$. Then $K_\Sigma$ is a compact open subgroup of $\G(\A_\fin)$ and $\Sigma=K_\Sigma\cap \G(\Q)$.
The map $\Sigma g\mapsto \G(\Q)g K_\Sigma$ induces a unitary $G$-isomorphism
$$
L^2(\Sigma\bs G)\tto\cong L^2(\G(\Q)\bs \G(\A)/K_\Sigma).
$$
Let $\G(\A)^+=\G(\A_\fin)\times G^+$.
This is a subgroup of index two in $\G(\A)$.
The above isomorphism restricts to a unitary $G^+$-isomorphism
$$
L^2(\Sigma\bs G^+)\tto\cong L^2(\G(\Q)\bs \G(\A)^+/K_\Sigma).
$$
Using the notation explained below, the space $S_2(\Sigma)$ can be identified with
$$
L_\cusp^2(\G(\Q)\bs \G(\A)^+/K_\Sigma)(\CD_2^+)(\eps_2).
$$
We explain the notation: 
The space of cusp forms in $L^2(\G(\Q)\bs\G(\A))$ is denoted by $L_\cusp^2(\G(\Q)\bs\G(\A))$.
The discrete series representations of $G^+$ are written as $\CD_2^+,\CD_2^-,\CD_4^+,\CD_4^-,\dots$ as in \cite{Knapp}.
Next, for $j\in\Z$ we denote by $\eps_{2j}$ the character of $K={\rm PSO}(2)$ given by $\eps_{2j}\smat{\cos\theta}{-\sin\theta}{\sin\theta}{\cos\theta}=e^{2ij\theta}$.
For any $K$-module $V$ we write $V(\eps_{2j})$ for its $\eps_{2j}$-isotypical component.
This explains the above notation.

\begin{lemma}\label{5.5}
Let $\Sigma\subset\Ga$ be a torsion-free congruence subgroup.
The order-lowering map
$$
\Gr S_{2,2}(\Sigma)\hookrightarrow\( S_2(\Sigma)\oplus \ol{S_2(\Sigma)}\)\otimes S_2(\Sigma)
$$
maps $\Gr S_{2,2}(\Sigma)$ to a subspace of codimension one.
\end{lemma}

\prf
The dimension of $S_2(\Sigma)$ is $g$, the genus of $\Sigma$, so the dimension of the right hand side is $2g^2$.
The dimension of $\Gr S_{2,2}(\Sigma)$ has been computed in \cite{DiamSim} to be equal to $2g^2-1$.
\qed

Taking limits, the lemma implies that $\psi$ has image of codimension at most one.

The space $\widehat S_2$ can be identified with $L_\cusp^2(\G(\Q)\bs\G(\A))(\CD_2^+)(\eps_2)$.
Likewise, $\widehat{\ol{S_2}}$ identifies with $L_\cusp^2(\G(\Q)\bs\G(\A))(\CD_2^-)(\eps_{-2})$.
Let $U=\(\widehat S_2\oplus\widehat{\ol{S_2}}\)\otimes \widehat S_2$. 
The closure $W\subset U$ of $\Im(\psi)$ is a closed subspace of codimension at most one.
We want to show $W=U$.
Assume this is not the case.
Then $\dim U/W=1$ and $U/W$ carries a unitary representation of $\G(\A_\fin)$.
The only finite dimensional irreducible unitary representation of $\G(\A_\fin)$ is the trivial representation.
This means that we have a non-zero, $\G(\A_\fin)$-equivariant  linear map $U\to\C$.
The space $L_\cusp^2(\G(\Q)\bs\G(\A))$ decomposes as a direct sum of irreducibles, say
$$
L_\cusp^2(\G(\Q)\bs\G(\A))=\bigoplus_{\pi\in\widehat{\G(\A))}} N_\cusp(\pi)\,\pi.
$$
Then, as a $\G(\A_\fin)$-representation,
$$
L_\cusp^2(\G(\Q)\bs\G(\A))(\CD_2^+)(\eps_2)=
\bigoplus_{\stack{\pi\in\widehat{\G(\A))}}{\pi_\infty=\CD_2^+}} N_\cusp(\pi)\,\pi_\fin.
$$
This implies that $U$ is a direct sum of tensor products of irreducible unitary representations.
Hence there exist irreducible $\G(\A_\fin)$-subrepresentations $\pi, \eta$ of $L_\cusp^2(\G(\Q)\bs\G(\A))(\CD_2^+)(\eps_2)$ or $L_\cusp^2(\G(\Q)\bs\G(\A))(\CD_2^-)(\eps_{-2})$, such that there is a non-zero continuous, $\G(\A_\fin)$-equivariant linear map 
$$
\Phi: \pi\hat\otimes\eta\ \to\ \C.
$$
Let $C$ be the kernel of $\Phi$.
Then the orthogonal complement, $C^\perp$, is one-dimensional, say equal to $\C T$.
The space $\pi\hat\otimes\eta$ can be viewed as the space of all Hilbert-Schmidt operators from $\pi^*$ to $\eta$ and $T\ne 0$ is such an operator which is $\G(\A_\fin)$ equivariant.
By Schur's Lemma it follows that $\eta\cong\pi^*$ and $T$ can be chosen unitary.
Since $T$ is also Hilbert-Schmidt, it follows that $\pi$ and $\eta$ are finite dimensional.
The only finite dimensional unitary representation of $\G(\A_\fin)$ is the trivial representation, hence $\pi=\eta=\triv$.
It follows that the representation $\triv\otimes\CD_2^+$ is automorphic, which is false.
So we have reached a contradiction, hence $U=W$ as claimed.
The theorem follows.
\qed

{\bf Remark.}
The current proof does not apply to $S_{2,q}$ for $q>2$, as, by \cite{DiamSim},
$$
\dim \Gr S_{2,q}(\Sigma)=\frac 12\((g+\sqrt{g^2-1})^{q}+(g-\sqrt{g^2-1})^{q}\),
$$
and so the codimension of $\Psi(\Gr S_{2,q}(\Sigma))$, which is $2^{q-1}g^{q}-\dim \Gr S_{2,q}(\Sigma)$, tends to infinity as the genus tends to infinity.
More precisely, as $g\to\infty$, this codimension grows like a constant times $g^{q}$.

{\bf Question.}
Does the order-lowering map generally give an isomorphism
$$
\widehat{\Gr S}_{k,q}\ \tto\cong\ \(\widehat S_2\oplus\widehat{\ol{S_2}}\)^{\otimes (q-1)}\otimes\widehat S_k ?
$$

We finish this section by showing that $\widehat{S_{k,q}}$ does not contain an irreducible subrepresentation.
This explains the fact that $L$-functions of higher order forms have no Euler product.
Using Proposition \ref{4.2} and Theorem \ref{4.3} the claim is implied by the following more general theorem.

\begin{theorem}
Let $S$ be an infinite set of primes and let $\A_S$ be the ring of $S$-adeles.
For any two nontrivial $\pi,\eta\in\widehat{\G(\A_S)}$ the representation $\pi\otimes\eta$ does not have an irreducible subrepresentation.
\end{theorem}

\prf
Assume $\pi\otimes\eta$ does contain an irreducible subrepresentation.
We first show that it must be infinite-dimensional.
Note that the only finite-dimensional unitary representation of $\G(\A_S)$ is the trivial representation, since the same holds for $\G(\Q_p)$ or every $p$.

\begin{proposition}
Let $\pi,\eta\in\widehat G$ for a topological group $G$.
If $\pi\otimes\eta$ contains the trivial representation, then $\pi$ and $\eta$ are finite-dimensional and dual to each other. 
\end{proposition}

\prf
Since the space of $\pi\otimes\eta$ is the space of Hilbert-Schmidt operators from $\pi^*$ to $\eta$, this is an easy consequence of Schur's Lemma.
\qed

So we know that $\pi\otimes\eta$ contains an infinite-dimensional irreducible subrepresentation.
Let $V_\pi$ and $V_\eta$ be the Hilbert spaces of $\pi$ and $\eta$.
According to the assumption there exists an infinite-dimensional irreducible unitary representation $(\tau,V_\tau)$ of $\G(\A_S)$ and a $G$-equivariant isometry $\psi: V_\tau\hookrightarrow V_\pi\hat\otimes V_\eta$.
Let $V_\pi'$ and $V_\eta'$ be the continuous duals.
We get a map $V_\pi'\otimes V_\eta'\otimes V_\tau\to \C$, given by
$$
\al\otimes\beta\otimes v\ \mapsto\ \al\otimes\beta(\psi(v)).
$$
This map is continuous and $G$-equivariant.
In order to prove the theorem, it suffices to show the following lemma.\

\begin{lemma}
Let $\pi,\eta,\tau\in\widehat{\G(\A_S)}\sm\{ 1\}$ and let 
$$
\chi: V_\pi\otimes V_\eta\otimes V_\tau\ \to\ C
$$
be a continuous, $G$-equivariant linear map.
Then $\chi=0$.
\end{lemma}

\prf
The representations $\pi$, $\eta$, $\tau$ are infinite tensor products
$$
\pi=\bigotimes_{p\in S}\pi_p,\quad
\eta=\bigotimes_{p\in S}\eta_p,\quad
\tau=\bigotimes_{p\in S}\tau_p,
$$
where $\pi_p, \eta_p, \tau_p\in\widehat{\G(\Q_p)}$.
Almost all of these are unramified.
Therefore, there exist a prime $p\in S$ such that all three $\pi_p, \eta_p, \tau_p$ are unramified and $\G(\Q_p)$ is not compact.
It suffices to show

\begin{lemma}
Let $p$ be a prime and let $\pi,\eta,\tau\in\widehat\G(\Q_p)\sm\{ 1\}$ be unramified.
Then every continuous, $G$-equivariant linear map
$$
\chi: V_\pi\otimes V_\eta\otimes V_\tau\ \to\ \C
$$
is zero.
\end{lemma}

Note that there can exist non-zero $G$-equivariant maps on the subspace $V_\pi^\infty\otimes V_\eta^\infty\otimes V_\tau^\infty$ of smooth vectors (see Theorem 1.2 of \cite{Prasad}).  These are, however, not continuous in the Hilbert-space topology.

\prf
For the length of this proof  we write $G=\G(\Q_p)$.
Let $A$ be the split torus of all diagonal matrices and let $N$ be the subgroup of all upper triangular matrices with ones on the diagonal.
Then $P=AN$ is a parabolic subgroup.
Let $V_\pi^\infty$ be the set of all smooth vectors in $V_\pi$.
Then $(\pi,V_\pi^\infty)$ is a smooth unramified representation of $G$.
Any such is parabolically induced by an unramified character of $A$.
More precisely, for $\mu\in\C$ and $a=\smat t\ \ 1$ we write $a^\mu=|t|^\mu$, where $|\cdot|$ is the $p$-adic absolute value.
Then there exists $\mu\in\C$ such that $V_\pi^\infty$ is isomorphic to the space of all smooth functions $f:G\to\C$ with $f(anx)=a^{\mu+1/2}f(x)$ on which $G$ acts via $\pi_\mu(y)f(x)=f(xy)$.

In order for $\pi_\mu$ to be unitarizable, there are two possibilities \cite{Bump}.
Either $\mu\in \R i$, then $\pi_\mu$ is called a \e{unitary principal series representation}, or $\mu\ne 0$ and $-1/2<\mu<1/2$, in which case $\pi_\mu$ is called a \e{complementary series representation}.
For our lemma, we are reduced to the case $(\pi,\eta,\tau)=(\pi_\mu,\pi_\nu,\pi_\la)$ for some $\mu,\nu,\la\in\C$.

The map $\chi$ induces a linear functional on the space of all continuous functions $\phi$ on $G\times G\times G$ such that
$$
\phi(a_1n_1x_1,a_2n_2x_2,a_3n_3x_3)=a_1^{\mu+1/2}a_2^{\nu+1/2}a_3^{\la+1/2}\phi(x_1,x_2,x_3).
$$
Since $G=PK$, this space of functions can be identified with the space $C((K\cap P\bs K)^3)$.
The functional $\chi$ is continuous on $L^2(M)$, where $M$ is the $K^3$-invariant measure on $(K\cap P\bs K)^3$.
The restriction of $\chi$ to the space of continuous functions is continuous with respect to the compact-open topology and thus is the integral of a complex-valued measure $m$ on $(K\cap P\bs K)^3$.

The $G$-space $(P\bs G)^3=(K\cap P\bs K)^3$ consists of 5 orbits:
\begin{itemize}
\item the closed orbit $(1,1,1)G\cong P\bs G$,
\item the intermediate orbits $(w,1,1)G, (1,w,1)G, (1,1,w)G\cong A\bs G$ where $w=\smat \ 11\ $, and
\item the open orbit $x_0G=(wn_0,w,1)G\cong G$, where $n_0\in N\sm\{ 1\}$ is arbitrary.
\end{itemize}
The space $V_\mu$ can be identified with $L^2(P\bs G)=L^2(K\cap P\bs K)$, the $L^2$-space taken with respect to a $K$-invariant measure.
In this identification, the $G$-action becomes
$$
\pi_\mu(g)\ph(k)=\ul a(kg)^{\mu+1/2}\ph(\ul k(kg)),
$$
where for $g\in G$ we write $g=\ul a(g)\ul n(g)\ul k(g)$ with $\ul a(g)\in A$, $\ul n(g)\in N$ and $\ul k(g)\in K$.
Note that this decomposition is only unique up to factors from $K\cap P$, but the expression above is independent of these ambiguities.
The space $C(P\bs G)$ of continuous functions on $P\bs G$ thus is a subspace of $V_\mu$ and so the space $C_c(x_0G)$ of compactly supported continuous functions on the open orbit $x_0G$ is a subspace of $V_\mu\otimes V_\nu\otimes V_\la\subset L^2((P\bs G)^3)$.
The restriction of $\chi$ to $C_c(x_0G)$ defines a $G$-invariant measure on the orbit $x_0G$.
We want to show that this measure is non-zero.
Assume the contrary.
Then $\chi$ defines an invariant measure concentrated on the remaining orbits.
These, however, constitute sets of measure zero with respect to $M$, so $\chi$ being continuous on $L^2(M)$ implies $m$ equals zero on the lower dimensional orbits.

Therefore there exists $c\in\C$ such that
$$
\chi(\phi)=c\int_G\phi(x_0g)\,dg.
$$
As $\chi$ is continuous on $L^2(M)$, the measure $m$ is absolutely continuous with respect to $M$, so by the Radon-Nikodym Theorem there exists a measurable function $w$ with $dm=wdM$.
The $L^2$-continuity of $\chi$ implies $w\in L^2(M)$, which trivially is equivalent to $|w|\in L^2(M)$.
For $\phi\in L^2(M)$ one has
$$
\chi(\phi)=\int_G \ul a(wn_0g)^{\mu+\frac12}\ul a(wg)^{\nu+\frac12}\ul a(g)^{\la+\frac12}\phi(\ul k(wn_0g),\ul k(wg),\ul k(g))\,dg.
$$
One sees that the absolute value $|w|$ does not depend on the parameters $\mu,\nu,\la$ as long as they are in $i\R$.
That means that if  for one triple one has $\chi\ne 0$ 
then this will be so for all triples and so one has a 
$G$-equivariant injection $\pi_{-\la}\hookrightarrow\pi_\mu\otimes\pi_\nu$ 
for every $\la\in i\R$ which contradicts the fact that 
$\pi_\mu\otimes\pi_\nu$ is a separable Hilbert space.
This cleans up the case when all parameters are imaginary.

Using the $ANK$-integration formula we get that $\chi(\phi)$ equals
$$
\int_{ANK}\ul a(wn_0an)^{\mu+\frac12}\ul a(wn)^{\nu+\frac12}a^{\la+\frac12}\ \phi(\ul k(wn_0an)k,\ul k(wn)k,k)\,da\,dn\,dk.
$$
Suppose now that $\la\in(-\frac12,\frac12)\sm\{ 0\}$.
As by Theorem 4.5.3 of \cite{Bump} we have $\pi_\la\cong\pi_{-\la}$, we can assume $\la>0$.
Let $\CO= x_0G$ be the open orbit and let $\CO^+$ be the set of all $x_0ank\in\CO$, where $a\in A$, $n\in N$, $k\in K$ such that $a=\smat y\ \ 1$ with $|y|>1$.
Let $\CO^-=\CO\sm\CO^+$ and let $w=w_\la$ be the Radon-Nikodym derivative mentioned above.
Then $\pi_\la$ injects into $\pi_\mu\otimes\pi_\nu$ if and only if $w_\la\in L^2(\CO,M)$.
The formula above shows that
$$
w_\la(x_0ank)=a^\la\tilde w(x_0ank), 
$$
where $\tilde w$ does not depend on $\la$.
Assume $\pi_\la\hookrightarrow \pi_\mu\otimes\pi_\nu$, then $w_\la\in L^2(\CO,M)$.
On the set $\CO^+$ we have $|w_{\la'}|\le|w_\la|$ for $\la'\le \la$, hence $w_{\la'}\in L^2(\CO^+,M)$ for $\la'\le \la$.
As $\pi_{-\la}\cong\pi_\la$, we get 
$w_{-\la}\in L^2(\CO,M)$ as well and so $w_\la'\in L^2(\CO^-,M)$ for 
$\la'\ge -\la$.
Taking things together it follows
$$
w_{\la'}\in L^2(\CO,M) \quad\text{for every}\quad -\la\le\la'\le\la.
$$
So it follows that $\pi_{\la'}$ is a subrepresentation of $\pi_{\mu}\otimes\pi_\nu$ for every $\la'\in [-\la,\la]$ which again contradicts the separability of $\pi_\mu\otimes\pi_\nu$.
The lemma and the theorem are proven.
\qed

\subsection{Hodge structure}
We keep considering the case $\CG=\SL_2(\Q)$ and $\Ga=\SL_2(\Z)$.
Lemma \ref{5.5} is complemented by the following  lemma.

\begin{lemma}\label{6.1}
For every congruence subgroup $\Sigma\subset\Ga$ the image of the order-lowering map $\Gr S_{2,2}(\Sigma)\hookrightarrow \(S_2(\Sigma)\oplus\ol{ S_2(\Sigma)}\)\otimes S_2(\Sigma)$ contains $S_2(\Sigma)\otimes S_2(\Sigma)$.
\end{lemma}

\prf
Let $f,g\in S_2(\Sigma)$.
Then $g$ defines a homomorphism $\chi_g:\Sigma\to\C$ by $\chi_g(\ga)=\int_z^{\ga z} g(w)\,dw$ which does not depend on the choice of $z\in\HH$.
We ave to show that there exists $F\in S_{2,1}(\Sigma)$ such that $F|_2(\ga -1)=\chi_g(\ga)f$.
For $z_0\in\HH$ set
$$
F(z)\df f(z)\int_{z_0}^z g(w)\,dw.
$$
Then
\begin{eqnarray*}
F|_2(\ga -1)(z) &=& f(z)\(\int_{z_0}^{\ga z}g(w)\,dw-\int_{z_0}^{z}g(w)\,dw\)\\
&=& f(z)\int_{z}^{\ga z}g(w)\,dw=\chi_g(\ga) f(z).
\end{eqnarray*}
\mathqed

There is a natural pure Hodge structure of weight 1 on $S_2(\Sigma)\oplus \ol{ S_2(\Sigma)}$ given by the fact that $S_2(\Sigma)$ can be identified with the space of holomorphic $1$-forms on the compactification $\Sigma\bs\HH\cup\{\text{cusps}\}$.

\begin{proposition}
There is a unique Hodge structure on the space $\Gr S_{2,2}(\Sigma)\oplus\ol{\Gr S_{2,2}(\Sigma)}$ which makes the order-lowering map 
$$
\Gr S_{2,2}(\Sigma)\oplus\ol{\Gr S_{2,2}(\Sigma)}\ \hookrightarrow\ \(S_2(\Sigma)\oplus\ol{S_2(\Sigma)}\)^{\otimes 2}
$$
a morphism of Hodge structures.
It is pure of weight 2.
\end{proposition}

\prf
Clear by Lemmas \ref{5.5} and \ref{6.1}.
\qed

\section{Higher order cohomology}
\subsection{Generalities}
Let now $R$ be an arbitrary commutative unital ring again.
For an $R$-module $M$ and a set $S$ we write $M^S$ for the $R$-module of all maps from $S$ to $M$.
Then $M^\emptyset$ is the trivial module $0$. 
Up to isomorphy, the module $M^S$ depends only on the cardinality of $S$.
It therefore makes sense to define $M^c$ for any cardinal number $c$ in this way.
Note that $J_q/J_{q+1}$ is a free $R$-module.
Define
$$
N_{\Ga,P}(q)\df \dim_R J_q/J_{q+1}.
$$
Then $N_{\Ga,P}(q)$ is a possibly infinite cardinal number.

\begin{lemma}\label{1.1}
\begin{enumerate}[\rm (a)]
\item For every $q\ge 1$ there is a natural exact sequence
\begin{multline*}
0\to H_{q}^0(\Ga,P,V)\to H_{q+1}^0(\Ga,P,V)\to H^0(\Ga,V)^{N_{\Ga,P}(q)}\to\\
\to H_{q}^1(\Ga,P,V)\to H_{q+1}^1(\Ga,P,V)\to H^1(\Ga,V)^{N_{\Ga,P}(q)}\to\dots\\
\dots\to H_{q}^p(\Ga,P,V)\to H_{q+1}^p(\Ga,P,V)\to H^p(\Ga,V)^{N_{\Ga,P}(q)}\to\dots
\end{multline*}
\item Suppose that for a given $p\ge 0$ one has $H^p(\Ga,V)=0$.
Then it follows $H_q^p(\Ga,P,V)=0$ for every $q\ge 1$.
In particular, if $V$ is acyclic as $\Ga$-module, then $H_q^p(\Ga,P,V)=0$ for all $p,q\ge 1$.
\end{enumerate}\end{lemma}

\prf
Consider the exact sequence
$$
0\to J_q/J_{q+1}\to A/J_{q+1}\to A/J_q\to 0.
$$
As an $A$-module, $J_q/J_{q+1}$ is isomorphic to a direct sum $\bigoplus_\al R_\al$ of copies of $R=A/I$.
So we conclude that for every $p\ge 0$,
$$
\Ext_A^p(J_q/J_{q+1},V)\ \cong\ \prod_\al\Ext_A^p(R,V)\ \cong\ H^p(\Ga,V)^{N_{\Ga,P}(q)}.
$$
The long exact $\Ext$-sequence induced by the above short sequence is
\begin{multline*}
0\to\Hom_A(A/J_q,V)\to\Hom_A(A/J_{q+1},V)\to \Hom_A(J_q/J_{q+1},V)\to\\
\to\Ext_A^1(A/J_q,V)\to\Ext_A^1(A/J_{q+1},V)\to \Ext_A^1(J_q/J_{q+1},V)\to\\
\to\Ext_A^2(A/J_q,V)\to\Ext_A^2(A/J_{q+1},V)\to \Ext_A^2(J_q/J_{q+1},V)\to\dots
\end{multline*}
This is the claim (a).
For (b) we proceed by induction on $q$.
For $q=1$ the claim follows from $H^p_1(\Ga,V)=H^p(\Ga,V)$.
Inductively, assume the claim proven for $q$ and $H^p(\Ga,V)=0$.
As part of the above exact sequence, we have the exactness of
$$
H_{q}^p(\Ga,P,V)\to H_{q+1}^p(\Ga,P,V)\to H^p(\Ga,V)^{N_{\Ga,P}(q)}.
$$
By assumption, we have $H^p(\Ga,V)^{N_{\Ga,P}(q)}=0$ and by induction hypothesis the module $H_{q}^p(\Ga,P,V)$ vanishes.
This implies $H_{q+1}^p(\Ga,P,V)=0$ as well.
\qed

\subsection{Sheaf cohomology}
Let $Y$ be a topological space which is path-connected and locally simply connected.
Let $C\to Y$ be a normal covering of $Y$.
Let $\Ga$ be the fundamental group of $Y$ and let $X\tto\pi Y$ be the universal covering.
The fundamental group $P$ of $C$ is a normal subgroup of $\Ga$.

For a sheaf $\CF$ on $Y$ define
$$
H_{q}^0(Y,C,\CF)\df H_{q}^0(\Ga,P,H^0(X,\pi^*\CF)),\qquad q\in\N.
$$
Let $\Mod(R)$ be the category of $R$-modules, let $\Mod_R(Y)$ be the category of sheaves of $R$-modules on $Y$, and let $\Mod_R(X)_\Ga$ be the category of sheaves over $X$ with an equivariant $\Ga$-action.
Then $H_{q}^0(Y,C,\cdot)$ is a left exact functor from $\Mod_R(Y)$ to $\Mod(R)$.
We denote its right derived functors by $H_{q}^p(Y,C,\cdot)$ for $p\ge 0$.

\begin{lemma} Assume that the universal cover $X$ is contractible.
\begin{enumerate}[\rm (a)]
\item For each $p\ge 0$ one has a natural isomorphism $H_{1}^{p}(Y,C,\CF)\ \cong\ H^p(Y,\CF)$.
\item If a sheaf $\CF$ is $H^0(Y,\cdot)$-acyclic, then it is $H_{q}^0(Y,C,\cdot)$-acyclic for every $q\in\N$.
\end{enumerate}
\end{lemma}

Note that part (b) allows one to use flabby or fine resolutions to compute higher order cohomology.

\prf
We decompose the functor $H^0_{q}(Y,C,\cdot)$ as 
$$
\Mod_R(Y)\tto{\pi^*}\Mod_R(X)_\Ga\tto{H^0(X,\cdot)}\Mod(R[\Ga])\tto{H^0_{q}(\Ga,P,\cdot)}\Mod(R).
$$
The functor $\pi^*$ is exact and maps injectives to injectives.
We claim that $H^0(X,\cdot)$ has the same properties.
For the exactness, consider the commutative diagram
$$
\xymatrix{
{\Mod_R(X)_\Ga}\ar[r]^{H^0}\ar[d]^f
	&{\Mod(R[\Ga])}\ar[d]^f\\
{\Mod_R(X)}\ar[r]^{H^0}
	&{\Mod(R),}
}
$$
where the vertical arrows are the forgetful functors.
As $X$ is contractible, the functor $H^0$ below is exact.
The forgetful functors have the property, that a sequence upstairs is exact if and only if its image downstairs is exact.
This implies that the above $H^0$ is exact.
It remains to show that $H^0$ maps injective objects to injective objects.
Let $\CJ\in\Mod_R(X)_\Ga$ be injective and consider a diagram with exact row in $\Mod(R[\Ga])$, 
$$
\xymatrix{
{0}\ar[r]
	&{M}\ar[r]\ar[d]^{\ph}
	&{N}\\
&{H^0(X,\CJ).}
}
$$
The morphism $\ph$ gives rise to a morphism $\phi:M\times X\to \CJ$, where $M\times X$ stands for the constant sheaf with stalk $M$.
Note that $H^0(X,\phi)=\ph$.
As $\CJ$ is injective, there exists a morphism $\psi:N\times X\to \CJ$ making the diagram
$$
\xymatrix{
{0}\ar[r]
	&{M\times X}\ar[r]\ar[d]^{\phi}
	&{N\times X}\ar[dl]^{\psi}\\
&{\CJ}
}
$$
commutative.
This diagram induces a corresponding diagram on the global sections, which implies that $H^0(X,\CJ)$ is indeed injective.

For a sheaf $\CF$ on $Y$ it follows that
$$
H^p(Y,\CF)=R^p(H^0(Y,\CF))=R^pH^0_{q}(\Ga,P,\CF)\circ H^0_\Ga\circ\pi^*= H_{q}^p(Y,C,\CF).
$$
Now let $\CF$ be acyclic.
Then we conclude $H_{0}^p(Y,C,\CF)=0$ for every $p\ge 1$, so the $\Ga$-module $V=H^0(X,\pi^*\CF)$ is $\Ga$-acyclic.
The claim follows from Lemma \ref{1.1}.
\qed

\subsection{Arithmetic groups}
We consider the case $R=\C$.
Let $G$ be a semisimple Lie group with compact center and let $X=G/K$ be its symmetric space.
Let $\Ga\subset G$ be an arithmetic subgroup which is torsion-free, and let $\Sigma\subset \Ga$ be a normal subgroup.
Let $Y=\Ga\bs X$, then $\Ga$ is the fundamental group of the manifold $Y$, and the universal covering $X$ of $Y$ is contractible.
This means that we can apply the results of the last section.

\begin{theorem}
Let $(\sigma,E)$ be a finite dimensional representation of $G$.
There is a natural isomorphism
$$
H_{q}^p(\Ga,P,E)\ \cong\ H_{\g,K}^p(H_{q}^0(\Ga,P, C^\infty(G))\otimes E),
$$
where the right hand side is the $(\g,K)$-cohomology, see \cite{BW}.
\end{theorem}

Note that, as a consequence of the definition of $(\g,K)$-cohomology and the fact that $E$ is finite dimensional, one can replace the module $H_{q}^0(\Ga,P,C^\infty(G))$ by its subspace of $K$-finite vectors
$$
H_{q}^0(\Ga,P,C^\infty(G))_K=H_{q}^0(\Ga,P,C^\infty(G)_K),
$$
where $K$ acts on $C^\infty(G)$ by right translations.

\prf
Let $\CF_E$ be the locally constant sheaf on $Y$ corresponding to $E$.
Let $\Omega_Y^p$ be the sheaf of complex valued $p$-differential forms on $Y$.
Then $\Omega_Y^p\otimes\CF_E$ is the sheaf of $\CF_E$-valued differential forms.
These form a fine resolution of $\CF_E$:
$$
0\to\CF_E\to\C^\infty\otimes\CF_E\tto{d\otimes 1}\Omega_Y^1\otimes\CF_E\to\dots
$$
Since $\pi^*\Omega_Y^\bullet=\Omega_X^\bullet$, 
we conclude that $H^p_{q}(\Ga,P,E)$ is the cohomology of the complex $H_{q}^0(\Ga,P,H^0(X,\Omega_X^\bullet\otimes E))$.
Let $\g$ and $\k$ be the Lie algebras of $G$ and $K$ respectively, and let $\g=\k\oplus\p$ be the Cartan decomposition.
Then $H^0(X,\Omega^p\otimes\CF_E)=(C^\infty(G)\otimes\bigwedge^p\p)^K\otimes E$. 
Mapping a form $\omega$ in this space to $(1\otimes x^{-1})\omega(x)$ one gets an isomorphism to $(C^\infty(G)\otimes\bigwedge^p\p\otimes E)^K$, where $K$ acts diagonally on all factors and $\Ga$ now acts on $C^\infty(G)$ alone.
The claim follows.
\qed

Let $U(\g)$ act on $C^\infty(G)$ as algebra of left invariant differential operators.
Let $\norm\cdot$ be a norm on $G$, see \cite{Wall}, Section 2.A.2.
Recall that a function $f\in C^\infty(G)$ is said to be \e{of moderate growth}, if for every $D\in U(\g)$ one has $Df(x)=O(\norm x^a)$ for some $a>0$.
The function $f$ is said to be of \e{uniform moderate growth}, if the exponent $a$ above can be chosen independent of $D$,
Let $C_\mg^\infty(G)$ and $C_\umg^\infty(G)$ denote the spaces of functions of moderate growth and uniform moderate growth respectively.

Let $\z$ be the center of the algebra $U(\g)$.
Let $\CA(G)$ denote the space of functions $f\in C^\infty(G)$ such that
\begin{itemize}
\item $f$ is of moderate growth,
\item $f$ is right $K$-finite, and
\item $f$ is $\z$-finite.
\end{itemize}

\begin{proposition}\label{2.2*}
For $\Omega=C_\umg^\infty(G),C_\mg^\infty(G),C^\infty(G)$ one has
$$
H_{q}^1(\Ga,P,\Omega)= 0
$$
for every $q\in\N$.
\end{proposition}

\prf
In order to prove the proposition, it suffices by Lemma \ref{1.1} (b), to consider the case $q=0$.
A 1-cocycle is a map $\al:\Ga\to \Omega$ such that $\al(\ga\tau)=\ga\al(\tau)+\al(\ga)$ holds for all $\ga,\tau\in\Ga$.
We have to show that for any given such map $\al$ there exists $f\in \Omega$ such that $\al(\tau)=\tau f-f$.
To this end 
consider the symmetric space $X=G/K$ of $G$.
Let $d(xK,yK)$ for $x,y\in G$ denote the distance in $X$ induced by the $G$-invariant Riemannian metric.
For $x\in G$ we also write $d(x)=d(xK,eK)$.
Then the functions $\log\norm x$ and $d(x)$ are equivalent in the sense that there exists a constant $C>1$ such that
$$
\frac 1C d(x)\ \le\ \log\norm x\ \le\ Cd(x)
$$
or
$$
\norm x\ \le\ e^{Cd(x)}\ \le\ \norm x^{C^2}
$$
holds for every $x\in G$.
We 
define
$$
\CF=\{ y\in G: d(y)< d(\ga y)\ \forall \ga\in\Ga\sm\{ e\}\}.
$$
As $\Ga$ is torsion-free, this is a fundamental domain for the $\Ga$ left translation action on $G$.
In other words, $\CF$ is open, its boundary is of measure zero, and there exists a set of representatives $R\subset G$ for the $\Ga$-action such that $\CF\subset R\subset \ol\CF$.
Next let $\ph\in C_c^\infty(G)$ with $\ph\ge 0$ and $\int_G\ph(x)\,dx=1$.
Set $u=\1_\CF *\ph$, where $\1_\CF$ is the indicator function of $\CF$ and $*$ is the convolution product $f*g(x)=\int_Gf(y)g(y^{-1}x)\,dy$.
Let $\CC$ be the support of $\ph$.
Then the support of $u$ is a subset of $\ol\CF\CC$ and the sum $\sum_{\tau\in\Ga} u(\tau^{-1} x)$ is locally finite in $x$.
For a function $h$ on $G$ and $x,y\in G$ we write $L_yh(x)=h(y^{-1} x)$.
Then for a convolution product one has $L_y(f*g)=(L_yf)*g$, and so
$$
\sum_{\tau\in\Ga}u(\tau^{-1}x)=\(\sum_{\tau\in\Ga}L_\tau\1_\CF\)*\ph.
$$
The sum in parenthesis is equal to one on the complement of a nullset.
Therefore,
$$
\sum_{\tau\in\Ga}u(\tau^{-1} x)\ \equiv\ 1.
$$
Set
$$
f(x)= -\sum_{\tau\in\Ga}\al(\tau)(x)\,u(\tau^{-1} x).
$$
\begin{lemma}
The function $f$ lies in the space $\Omega$.
\end{lemma}

\prf
Since the sum is locally finite, 
the function $f$ is smooth and for each $D\in U(\g)$ we have $Df=D\sum_\tau \al(\tau)L_\tau u=\sum_\tau D(\al(\tau)L_\tau u)$.
By the Leibniz rule, $D(\al(\tau)L_\tau u)$ is a finite sum of terms of the form $(D\al(\tau))(D'L_\tau u)$ with $D,D'\in U(\g)$.
We have to show that the sum
$$
h=\sum_{\tau\in\Ga}(D\al(\tau))(D'L_\tau u)
$$
is in $\Omega$.
We will treat the case $\Omega= C_\umg^\infty(G)$ here, the case $C_\mg^\infty$ is similar and the case $C^\infty(G)$ is trivial, as no growth bounds are required.
We first observe that
$$
D'(L_\tau u)=(L_\tau \1_\CF)*(D'\ph).
$$
Set
$$
S=\{ \ga\in\Ga\sm\{ e\} : \ga\ol\CF\cap\ol\CF\ne \emptyset\}.
$$
Then $S$ is a finite symmetric generating set for $\Ga$.
For $\ga\in\Ga$, let $\CF_\ga$ be the set of all $x\in G$ with $d(x) <d(\ga x)$.
Then
$$
\CF=\bigcap_{\ga\in\Ga\sm\{ e\}}\CF_\ga
$$
Let $\tilde\CF=\bigcap_{s\in S}\CF_s$. We claim that $\CF=\tilde \CF$.
As the intersection runs over fewer elements, one has $\CF\subset\tilde\CF$.
For the converse note that for every $s\in S$ the set $s\ol\CF/K$ lies in $X\sm\tilde\CF/K$, therefore $\CF/K$ is a connected component of $\tilde\CF/K$.
By the invariance of the metric, we conclude that $x\in\CF_\ga$ if and only if $d(xK,eK)<d(xK,\ga^{-1}K)$.
This implies that $\CF_\ga/K$ is a convex subset of $X$.
Any intersection of convex sets remains convex, therefore $\tilde\CF/K$ is convex and hence connected, and so $\tilde \CF/K=\CF/K$, which means $\tilde\CF=\CF$.

Likewise we get $\ol\CF=\bigcap_{s\in S}\ol{\CF_s}$.
The latter implies that for each $x\in G\sm\ol\CF$ there exists $s\in S$ such that $d(s^{-1}x)<d(x)$.
Iterating this and using the fact that the set of all $d(\ga x)$ for $\ga\in\Ga$ is discrete, we find for each $x\in G\sm\ol\CF$ a chain of elements $s_1,\dots,s_n\in S$ 
such that $d(x) >d(s_1^{-1}x)>\dots>d(s_n^{-1}\cdots s_1^{-1} x)$
 and $s_n^{-1}\cdots s_1^{-1} x\in\ol\CF$.
The latter can be written as $x\in s_1\dots s_n\ol\CF$.
Now let $\tau\in\Ga$ and suppose $u(\tau^{-1} x)\ne 0$.
Then $x\in\ol\CF C$, so, choosing $C$ small enough, we can assume $x\in s\tau\ol\CF$ for some $s\in S\cap\{ e\}$.
As the other case is similar, we can assume $s=e$.
It suffices to assume $x\in\tau\CF$, as we only need the estimates on the dense open set $\Ga\CF$.
So then it follows $\tau= s_1\dots s_n$.

Let $D\in U(\g)$.
As $\al$ maps to $\Omega= C_\umg^\infty(G)$, for every $\ga\in\Ga$ there exist $C(D,\ga), a(\ga)>0$ such that
$$
|D\al(\ga)(x)|\ \le\ C(D,\ga)\,\norm x^{a(\ga)}.
$$
The cocycle relation of $\al$ implies
$$
\al(\tau)(x)=\sum_{j=1}^n\al(\ga_j)(s_{j-1}^{-1}\dots s_1^{-1} x).
$$
We get
\begin{eqnarray*}
|D\al(\tau)(x)| &\le& \sum_{j=1}^n C(D,s_j)\norm{s_{j-1}^{-1}\dots s_1^{-1}x}^{a(s_j)}\\
&\le& \sum_{j=1}^n C(D,s_j)\, e^{Cd(s_{j-1}^{-1}\dots s_1^{-1}x)a(s_j)}\\
&\le& \sum_{j=1}^n C(D,s_j)\, e^{Cd(x)a(s_j)}\\
&\le& \sum_{j=1}^n C(D,s_j)\norm{x}^{C^2a(s_j)}\ \le\ nC_0(D)\norm x^{a_0},
\end{eqnarray*}
where $C(D)=\max_j C(D,s_j)$ and $a_0=C^2\max_jd(s_j)$.
It remains to show that $n$ only grows like a power of $\norm x$.
To this end let for $r>0$ denote $N(r)$ the number of $\ga\in\Ga$ with $d(\ga)\le r$.
Then a simple geometric argument shows that
$$
N(r)=\frac 1{\vol\CF}\vol\(\bigcup_{\ga:d(\ga)\le r}\ga\CF/K\)\ \le\ C_1\vol (B_{2r}),
$$
where $B_{2r}$ is the ball of radius $2r$ around $eK$.
Note that for the homogeneous space $X$ there exists a constant $C_2>0$ such that $\vol B_{2r}\le e^{C_2 r}$.
Now $n\le N(d(x))$ and therefore
$$
n\ \le\ C_1\vol B_{2d(x)}\ \le\ C_1e^{C_2 d(x)}\ \le\ C_1 \norm x^{C_3}
$$
for some $C_3>0$.
Together it follows that there exists $C(D)>0$ and $a>0$ such that
$$
|D\al(\tau)(x)|\ \le\ C(D)\,\norm x^a.
$$
Therefore it follows that 
\begin{eqnarray*}
|h(x)| &\le& \sum_{\tau\in\Ga} |D\al(\tau)(x)||L_\tau\1_\CF*D'\ph(x)|\\
&\le& C(D)\norm x^a\sum_{\tau\in\Ga}|L_\tau\1_\CF*D'\ph(x)|\\
&\le& C(D)\norm x^a\sum_{\tau\in\Ga}L_\tau\1_\CF*|D'\ph|(x)\\
&=& C(D)\norm x^a\(\sum_{\tau\in\Ga}L_\tau\1_\CF\)*|D'\ph|(x)\\
&=& C(D)\norm x^a\underbrace{\1*|D'\ph|(x)}_{=\rm const.}\\
\end{eqnarray*}
This is the desired estimate which shows that $f\in\Omega$.
The lemma is proven.
\qed

To finish the proof of the proposition, we now compute for $\ga\in\Ga$,
\begin{eqnarray*}
\ga f(x)- f(x) &=& f(\ga^{-1} x)-f(x)\\
&=& \sum_{\tau\in\Ga}\al(\tau x) u(\tau^{-1} x) -
\al(\tau)(\ga^{-1} x)u(\tau^{-1}\ga^{-1} x)\\
&=& \sum_{\tau\in\Ga}\al(\tau) (x) u(\tau^{-1}x) 
+\al(\ga)(x)\sum_{\tau\in\Ga} u((\ga\tau)^{-1}x)\\
&& -\sum_{\tau\in\Ga}\al(\ga \tau)(x)u((\ga\tau)^{-1} x)
\end{eqnarray*}
The first and the last sum cancel and the middle sum is $\al(\ga)(x)$.
Therefore, the  proposition is proven.
\qed

\begin{proposition}\label{2.4*}
For every $q\ge 1$ there is an exact sequence of continuous $G$-homomorphisms,
$$
0\to H_{q}^0(\Ga,P,C_*^\infty(G))\tto\phi  H_{q+1}^0(\Ga,P,C_*^\infty(G))
\tto\psi C_*^\infty(\Ga\bs G)^{N_{\Ga,P}(q)}\to 0,
$$
where $\phi$ is the inclusion map and $*$ can be $\emptyset, \umg$, or $\mg$.
\end{proposition}

\prf
This follows from Lemma \ref{1.1} together with Proposition \ref{2.2*}.
\qed

The space $C^\infty(G)$ carries a natural topology which makes it a nuclear topological vector space.
For every $q\in\N$, the space $H_{q}^0(\Ga,P,C^\infty(G))$ is a closed subspace.
If $\Ga$ is cocompact, then one has the isotypical decomposition
$$
H_{1}^0(\Ga,P, C^\infty(G))= C^\infty(\Ga\bs G)= \ol{\bigoplus_{\pi\in\hat G}} C^\infty(\Ga\bs G)(\pi),
$$
and $C^\infty(\Ga\bs G)(\pi)\cong m_\Ga(\pi)\pi^\infty$,
where the sum runs over the unitary dual $\hat G$ of $G$, and for $\pi\in\hat G$ we write $\pi^\infty$ for the space of smooth vectors in $\pi$.
The multiplicity $m_\Ga(\pi)\in\N_0$ is the multiplicity of $\pi$ as a subrepresentation of $L^2(\Ga\bs G)$, i.e.,
$$
m_\Ga(\pi)=\dim\Hom_G\(\pi,L^2(\Ga\bs G)\).
$$
Finally, the direct sum $\ol{\bigoplus}$ means the closure of the algebraic direct sum in $C^\infty(G)$.
We write $\hat G_\Ga$ for the set of all $\pi\in\hat G$ with $m_\Ga(\pi)\ne 0$.

Let $\pi\in\hat G$.
A smooth representation $(\beta,V_\beta)$ of $G$ is said to be \e{of type $\pi^\infty$}, if it is of finite length and every irreducible subquotient is isomorphic to $\pi^\infty$.
For a smooth representation $(\eta,V_\eta)$ we define the \e{$\pi^\infty$-isotype} as
$$
V_\eta(\pi^\infty)\df \ol{\sum_{\stack{V_\beta\subset V_\eta}{\beta{\rm\ of\ type\ }\pi}}}V_\beta,
$$
where the sum runs over all subrepresentations $V_\beta$ of type $\pi^\infty$.

\begin{theorem}
Suppose $\Ga$ is cocompact and let $*\in\{ \emptyset,\mg,\umg\}$.
We write $V_q=H_{q}^0(\Ga,P,C_*^\infty(G))$.
For every $q\ge 1$ there is an isotypical decomposition
$$
V_q=\ol{\bigoplus_{\pi\in\hat G_\Ga}}V_q(\pi),
$$
and each $V_q(\pi)$ is of type $\pi$ itself.
The exact sequence of Proposition \ref{2.4*} induces an exact sequence
$$
0\to V_{q}(\pi)\to V_{q+1}(\pi)\to (\pi^\infty)^{m_\Ga(\pi)N_{\Ga,P}(q)}\to 0
$$
for every $\pi\in\hat G_\Ga$.
\end{theorem}

\prf
We will prove the theorem by reducing to a finite dimensional situation by means of considering infinitesimal characters and $K$-types.
For this let $\hat\z=\Hom(\z,\C)$ be the set of all algebra homomorphisms from $\z$ to $\C$.
For a $\z$-module $V$ and $\chi\in\hat\z$ let
$$
V(\chi)\df \{ v\in V: \forall z\in\z\ \exists n\in\N\ (z-\chi(z))^nv=0\}
$$
be the \e{generalized $\chi$-eigenspace}.
Since $\z$ is finitely generated, one has
$$
V(\chi)= \{ v\in V: \exists n\in\N\ \forall z\in\z\ (z-\chi(z))^nv=0\}.
$$
For $\chi\ne \chi'$ in $\z$ one has $V(\chi)\cap V(\chi')=0$.
Recall that the algebra $\z$ is free in $r$ generators, where $r$ is the absolute rank of $G$.
Fix a set of generators $z_1,\dots, z_r$.
The map $\chi\mapsto (\chi(z_1),\dots,\chi(z_r))$ is a bijection $\hat\z\to\C^r$.
We equip $\hat\z$ with the topology of $\C^r$.
This topology does not depend on the choice of the generators $z_1,\dots,z_r$.

Let $\Ga\subset G$ be a discrete cocompact subgroup. 
Let $\hat\z_\Ga$ be the set of all $\chi\in\hat\z$ such that the generalized eigenspace $C^\infty(\Ga\bs G)(\chi)$ is non-zero.
The $\hat\z_\Ga$ is discrete in $\hat\z$, more sharply there exists 
$\eps_\Ga>0$ such that for any two $\chi\ne\chi'$ in $\hat\z_\Ga$ 
there is $j\in\{ 1,\dots,r\}$ such that $|\chi(z_j)-\chi'(z_j)|>\eps_\Ga$.

\begin{proposition}
Let $*\in\{ \emptyset,\mg,\umg\}$.
For every $q\ge 0$ and every $\chi\in\hat\z$ the space 
$V_q(\chi)=H_{q,\Sigma}^0(\Ga,C^\infty_*(G))(\chi)$ coincides with
$$
\bigcap_{z\in\z}\ker(z-\chi(z))^{2^{q-1}},
$$
and is therefore a closed subspace of $V_q$.
The representation of $G$ on $V_q(\chi)$ is of finite length.

The space $V_q(\chi)$ is non-zero only if $\chi\in\hat\z_\Ga$.
One has a decomposition
$$
H_{q}^0(\Ga,P,C^\infty_*(G))=\ol{\bigoplus_{\chi\in\hat\z_\Ga}}H_{q}^0(\Ga,P,C^\infty_*(G))(\chi).
$$
The exact sequence of Proposition \ref{2.4*} induces an exact sequence
$$
0\to V_{q}(\chi)\to  V_{q+1}(\chi)
\to \bigoplus_{\pi\in\hat G_\chi}m_\Ga(\pi)N_{\Ga,P}(q)\pi
\to 0.
$$
\end{proposition}

\prf
All assertions are clear for $q=1$.
We proceed by induction.
Fix $\chi\in\hat\z_\Ga$.
The  exact sequence
$$
0\to V_{q}\to V_{q+1}\tto{\psi}V_1^{N_{\Ga,P}(q)}\to 0.
$$
induces an exact sequence
$$
0\to V_{q}(\chi)\to V_{q+1}(\chi)\tto{\psi_\chi}V_1(\chi)^{N_{\Ga,P}(q)}.
$$
Let $v\in V_1(\chi)^{N_{\Ga,P}(q)}$.
As $\psi$ is surjective, one finds $u\in V_{q+1}$ with $\psi(u)=v$.
We have to show that one can choose $u$ to lie in $V_{q+1}(\chi)$.
We have $(z-\chi(z))v=0$ for every $z\in\z$.
Therefore $(z-\chi(z))u\in V_{q}$.
Inductively we assume the decomposition to holds for $V_{q}$, so we can write
$$
(z_j-\chi(z_j))u=\sum_{\chi'\in\hat\z_\Ga}u_{j,\chi'},
$$
for $1\le j\le r$ and $u_{j,\chi'}\in\ker(z-\chi'(z))^{2^{q-1}}$ for every $z\in\z$.
For every $\chi'\in\hat\z_\Ga\sm\{\chi\}$ we fix some index $1\le j(\chi')\le r$ with $|\chi(z_{j(\chi')})-\chi'(z_{j(\chi')})|>\eps_\Ga$.
On the space
$$
\ol{\bigoplus_{\chi':j(\chi')=j}}V_q(\chi')
$$
the operator $z_j-\chi(z_j)$ is invertible and the inverse $(z_j-\chi(z_j))^{-1}$ is continuous.
We can replace $u$ with
$$
u-\sum_{\chi'\in\hat\z_{\Ga}\sm\{\chi\}}(z_{j(\chi')}-\chi(z_{j(\chi')}))^{-1}u_{j(\chi'),\chi'}.
$$
We end up with $u$ satisfying $\psi(u)=v$ and
$$
(z_1-\chi(z_1))\cdots (z_r-\chi(z_r))u\ \in\ V_{q}(\chi)=\bigcap_{z\in\z}\ker(z-\chi(z))^{2^{q-1}}.
$$
So for every $z\in\z$ one has
$$
0= (z_1-\chi(z_1))\cdots (z_r-\chi(z_r))(z-\chi(z))^{2^{q-1}}u,
$$
which implies
$$
(z-\chi(z))^{2^{q-1}}u\in\ker\((z_1-\chi(z_1))\cdots (z_r-\chi(z_r))\).
$$
As the set $\z_\Ga$ is countable, one can, depending on $\chi$, choose the generators $z_1,\dots,z_r$ in a way that $\chi(z_j)\ne \chi'(z_j)$ holds for every $j$ and every $\chi'\in\hat\z_\Ga\sm\{\chi\}$.
Therefore the operator $(z_1-\chi(z_1))\cdots (z_r-\chi(z_r))$ 
is invertible on $V_{q}(\chi')$ for every 
$\chi'\in\hat\z_\Ga\sm\{\chi\}$ and it follows 
$(z-\chi(z))^{2^{q-1}}u\in V_{q}(\chi)\subset \ker(z-\chi(z))^{2^{q-1}}$ 
and therefore $u\in\ker((z-\chi(z))^{2^q})$. 
Since this holds for every $z$ it follows $u\in V_{q+1}(\chi)$ and hence $\psi_\chi$ is indeed surjective. One has an exact sequence
$$
0\to V_{q}(\chi)\to V_{q+1}(\chi)\to V_1(\chi)^{N_{\Ga,P}(q)}\to 0.
$$
Taking the sum over all $\chi\in\z_\Ga$ we arrive at an exact sequence
$$
0\to V_{q}\to \ol{\bigoplus_{\chi\in\hat\z_\Ga}}V_{q+1}(\chi)\to V_1^{N_{\Ga,P}(q)}\to 0.
$$
Hence we get a commutative diagram with exact rows
$$
\xymatrix{
0\ar[r] &V_q\ar[r]\ar[d]^=&
\ds\ol{\bigoplus_{\chi\in\hat\z_\Ga}}V_{q+1}(\chi)\ar[r]\ar@{^{(}->}_i[d]
&V_1^{N_{\Ga,P}(q)}\ar[r]\ar[d]^=&0\\
0\ar[r] & V_q\ar[r] & V_{q+1}\ar[r]& V_1^{N_{\Ga,P}(q)}\ar[r]& 0,
}
$$
where $i$ is the inclusion.
By the 5-Lemma, $i$ must be a bijection.
The proposition follows.
\qed

We now finish the proof of the theorem.
We keep the notation $V_q$ for the space $\H_q^0(\Ga,\Sigma,C_*^\infty(G))$.
For a given $\chi\in \hat\z_\Ga$ the $G$-representation $V_q(\chi)$ is of finite length, so the $K$-isotypical decomposition
$$
V_q(\chi)=\ol{\bigoplus_{\tau\in\hat K}}V_q(\chi)(\tau)
$$
has finite dimensional isotypes, i.e., $\dim V_q(\chi)(\tau)<\infty$.
Let $U(\g)^K$ be the algebra of all $D\in U(\g)$ such that $\Ad(k)D=D$ for every $k\in K$.
Then the action of $D\in U(\g)$ commutes with the action of each $k\in K$, and so $K\times U(\g)^K$ acts on every smooth $G$-module.
For $\pi\in\hat G$ the $K\times U(\g)^K$-module $V_\pi(\tau)$ is irreducible and $V_\pi(\tau)\cong V_{\pi'}(\tau')$ as a $K\times U(\g)^K$-module implies 
$\pi=\pi'$ and $\tau=\tau'$, see \cite{Wall}, Proposition 3.5.4.
As $V_q(\chi)(\tau)$ is finite dimensional. one gets
$$
V_q(\chi)(\tau)=\bigoplus_{\stack{\pi\in\hat G}{\chi_\pi=\chi}}V_q(\chi)(\tau)(\pi),
$$
where $V_q(\chi)(\tau)(\pi)$ is the largest $K\times U(\g)^K$-submodule of $V_q(\chi)(\tau)$ with the property that every irreducible subquotient is isomorphic to $V_\pi(\tau)$.
Let 
$$
V_q(\pi)=\ol{\bigoplus_{\tau\in\hat K}}V_q(\chi_\pi)(\tau)(\pi).
$$
The claims of the theorem follow from the proposition.
\qed

\subsection{The higher order Borel conjecture}

Let $(\sigma,E)$ be a finite dimensional representation of $G$. In \cite{Borel}, 
A. Borel has shown that the inclusions $C_\umg^\infty(G)\hookrightarrow C_\mg^\infty(G)\hookrightarrow C^\infty(G)$ induce isomorphisms in cohomology:
\begin{multline*}
H_{\g,K}^p(H^0(\Ga,C_\umg^\infty(G))\otimes E)\tto\cong 
H_{\g,K}^p(H^0(\Ga,C_\mg^\infty(G))\otimes E)\\
\tto\cong 
H_{\g,K}^p(H^0(\Ga,C^\infty(G))\otimes E).
\end{multline*}

In \cite{Franke}, J. Franke proved a conjecture of Borel stating that the inclusion $\CA(G)\hookrightarrow C^\infty(G)$ induces an isomorphism
$$
H_{\g,K}^p(H^0(\Ga,\CA(G))\otimes E)\tto\cong 
H_{\g,K}^p(H^0(\Ga,C^\infty(G))\otimes E).
$$

\begin{conjecture}
[Higher order Borel conjecture]
For $q\in\N$, the inclusion $\CA(G)\hookrightarrow C^\infty(G)$ induces an isomorphism
$$
H_{\g,K}^p(H_{q}^0(\Ga,P,\CA(G))\otimes E)\tto\cong 
H_{\g,K}^p(H_{q}^0(\Ga,P,C^\infty(G))\otimes E).
$$
\end{conjecture}

We can prove the higher order version of Borel's result.

\begin{theorem}
For each $q\in\N$, the inclusions $C_\umg^\infty(G)\hookrightarrow C_\mg^\infty(G)\hookrightarrow C^\infty(G)$ induce isomorphisms in cohomology:
\begin{multline*}
H_{\g,K}^p(H_{q}^0(\Ga,P,C_\umg^\infty(G))\otimes E)\tto\cong 
H_{\g,K}^p(H_{q}^0(\Ga,P,C_\mg^\infty(G))\otimes E)\\
\tto\cong 
H_{\g,K}^p(H_{q}^0(\Ga,P,C^\infty(G))\otimes E).
\end{multline*}
\end{theorem}

\prf
Let $\Omega$ be one of the spaces $\C_\umg^\infty(G)$ or $\C_\mg^\infty(G)$.

By Proposition \ref{2.4*} we get an exact sequence
$$
0\to H_{q}^0(\Ga,P,\Omega)\to  H_{q+1}^0(\Ga,P,\Omega)
\to H^0(\Ga,\Omega)^{N_{\Ga,P}(q)}\to 0,
$$
and the corresponding long exact sequences in $(\g,K)$-cohomology.
For each $p\ge 0$ we get a commutative diagram with exact rows
{\scriptsize
$$
\xymatrix{
{H_{\g,K}^p(H_{q}^0(\Ga,P,\Omega)\otimes E)}\ar[r]\ar[d]^{\al}
	&{H_{\g,K}^p(H_{q+1}^0(\Ga,P,\Omega)\otimes E)}\ar[r]\ar[d]^{\beta}
	&{H_{\g,K}^p(H^0(\Ga,\Omega)\otimes E)^{N_{\Ga,P}(q)}}\ar[d]^{\ga}\\
{H_{\g,K}^p(H_{q}^0(\Ga,P,C^\infty(G))\otimes E)}\ar[r]
	&{H_{\g,K}^p(H_{q+1}^0(\Ga,P,C^\infty(G))\otimes E)}\ar[r]
	&{H_{\g,K}^p(H^0(\Ga,C^\infty(G))\otimes E)^{N_{\Ga,P}(q)}.}
}
$$}
Borel has shown that $\ga$ is an isomorphism and that $\al$ is an isomorphism for $q=0$.
We prove that $\beta$ is an isomorphism by induction on $q$.
For the induction step we can assume that $\al$ is an isomorphism.
Since the diagram continues to the left and right with copies of itself where $p$ is replaced by $p-1$ or $p+1$, we can deduce that $\beta$ is an isomorphism by the 5-Lemma.
\qed

{\small Mathematisches Institut\\
Auf der Morgenstelle 10\\
72076 T\"ubingen\\
Germany\\
\tt deitmar@uni-tuebingen.de}

\today

\end{document}